\documentclass{amsart}
\usepackage{graphicx}

\newtheorem{theorem}{Theorem}[section]
\newtheorem{lemma}[theorem]{Lemma}
\newtheorem{corollary}[theorem]{Corollary}

\theoremstyle{definition}

\newtheorem{claim}[theorem]{Claim}
\newtheorem{conjecture}[theorem]{Conjecture}

\theoremstyle{remark}

\numberwithin{equation}{section}

\begin{document}

\title{Non-triviality of generalized alternating knots}

\author{Makoto Ozawa}
\address{Natural Science Faculty, Faculty of Letters, Komazawa University,
 1-23-1 Komazawa, Setagaya-ku, Tokyo, 154-8525, Japan}
\email{w3c@komazawa-u.ac.jp}





\begin{abstract}
In this article, we consider alternating knots on a closed surface in the 3-sphere, and show that these are not parallel to any closed surface disjoint from the prescribed one.
\end{abstract}

\maketitle

\section{Introduction}
How do you knit a rope to have a knot so that it does not come loose?
It is the simplest method to make the knot alternating, that is, a knot which possesses a knot diagram in which crossings alternate between under- and over-passes.

Since a long time before, alternating knots have been known to be non-trivial.
It is a monument in knot theory that a reducible, alternating knot diagram gives a minimal genus Seifert surface (\cite{C}, \cite{Mu1}) and a minimal crossing number (\cite{K}, \cite{Mu2}, \cite{T}).
By using a spanning disk, Menasco and Thistlethwaite directly proved that alternating knots are non-trivial (\cite{MT}).

As another way to prove the non-triviality of knots, it is effective to construct incompressible surfaces in the knot exterior.
To show that the asphericity of alternating link complements, Aumann proved that an alternating link can be contained in a closed surface obtained from a checkerboard coloring so that its complement is incompressible (\cite{A}).
The incompressibility of the checkerboard surface for an alternating knot was also used to prove Property P (\cite{DR}).

It is a purpose of this article to generalize the Aumann's result and to show the non-triviality of generalized alternating knots.

\section{Result}

Let $K$ be a knot in the 3-sphere $S^3$ and $S$ a closed surface embedded in $S^3-K$.
We say that $K$ is {\em trivial} for $S$ if it is ambient isotopic to a knot contained in $S$, and that $K$ is {\em parallel} to $S$ if there exists an annulus $A$ embedded in $S^3$ such that $A\cap K=\partial A\cap K=K$ and $A\cap S=\partial A-K$.
It follows that if $K$ is parallel to $S$, then $K$ is trivial for $S$, but the converse does not hold generally.
For example, consider the trefoil knot $K$ contained in a 3-ball $B$, and let $V$ be an unknotted solid torus containing $B$.
Then, $K$ is trivial for $\partial V$, but $K$ is not parallel to $\partial V$.
In the case that $S$ is a 2-sphere, it holds that $K$ is trivial for $S$ if and only if $K$ is parallel to $S$.

Let $F$ be a closed surface of positive genus in $S^3$ and $K$ a knot contained in $F$.
The {\em representativity} $r(F,K)$ of a pair $(F,K)$ is defined as the minimal number of intersecting points of $K$ and $\partial D$, where $D$ ranges over all compressing disks for $F$ in $S^3$.
We note that there exists at least one compressing disk for $F$ in $S^3$ since $S^3$ does not contain closed incompressible surfaces.
This definition is a modification of the representativity for a graph embedded in a surface (\cite{RV}).
It follows from Lemma \ref{representativity} that $r(F,K)\ge1$ if and only if $F\cap E(K)$ is incompressible in $E(K)$, and $r(F,K)\ge2$ if and only if $F\cap E(K)$ is incompressible and $\partial$-incompressible in $E(K)$, where $E(K)$ denotes the exterior of $K$ in $S^3$.

The following theorem is a mild generalization of Corollary 2 in \cite{IO}.

\begin{theorem}\label{non-parallel}
Let $K$ be a knot contained in a closed surface $F$ of positive genus in $S^3$ and $S$ be a closed surface disjoint from $F$.
If $r(F,K)\ge 2$ and $F$ is compressible in $S^3-S$ on both sides, then $K$ is not parallel to $S$.
\end{theorem}

Let $F$ be a closed surface embedded in $S^3$ and $K$ a knot contained in $F\times [-1,1]$.
Suppose that $\pi(K)$ is a regular projection on $F$, where $\pi : F\times [-1,1]\to F\times \{0\}=F$ is the projection.
Then, we have a regular diagram on $F$ obtained from $\pi(K)$ by adding the over/under information to each double point, and we denote it by the same symbol $\pi(K)$ in this article.
As usual, a diagram $\pi(K)$ on $F$ is said to be {\em alternating} if it has alternating over- and under-crossings as the diagram $\pi(K)$ is traversed on $F$.
We say that a diagram $\pi(K)$ on $F$ is {\em reduced} if there is no disk region of $F-\pi(K)$ which meets only one crossing.
We say that a diagram $\pi(K)$ on $F$ is {\em prime} if it contains at least one crossing and for any loop $l$ intersecting $\pi(K)$ in two points except for crossings, there exists a disk $D$ in $F$ such that $D\cap \pi(K)$ consists of an embedded arc.

The following theorem was announced as Theorem 5 in \cite{IO}.

\begin{theorem}\label{alternating}
Let $F$ be a closed surface embedded in $S^3$, $K$ a knot contained in $F\times [-1,1]$ which has a reduced, prime, alternating diagram on $F$.
Then, we have the following.
\begin{enumerate}
	\item $F-\pi(K)$ consists of open disks.
	\item $F-\pi(K)$ admits a checkerboard coloring.
	\item $K$ bounds a non-orientable surface $N$ coming from the checkerboard coloring.
	\item $K$ can be isotoped into $\partial N(N)$ so that $\partial N(N)-K$ is connected.
	\item $r(\partial N(N),K)\ge 2$.
\end{enumerate}
\end{theorem}

Theorem \ref{alternating} assures the existence of an incompressible and $\partial$-incompressible separating orientable surface of integral boundary slope in the exterior of a generalized alternating knot.
The following corollary slightly extends Corollary 4 in \cite{IO}.

\begin{corollary}\label{alternating2}
Let $F$ be a closed surface embedded in $S^3$, $K$ a knot contained in $F\times [-1,1]$ which has a reduced, prime, alternating diagram on $F$.
Then, $K$ is not parallel to any closed surface in $S^3-(F\times[-1,1])$, and $K$ is not trivial in particular.
\end{corollary}

We can also consider the link case.

\begin{theorem}\label{alternating link}
Let $F$ be a closed surface embedded in $S^3$, $L$ a link contained in $F\times [-1,1]$ which has a reduced, prime, alternating diagram on $F$.
Then, we have the following.
\begin{enumerate}
	\item $F-\pi(L)$ consists of open disks.
	\item $F-\pi(L)$ admits a checkerboard coloring.
	\item $L$ bounds a connected surface $N$ coming from the checkerboard coloring.
	\item $L$ can be isotoped into $\partial N(N)$.
	\item $r(\partial N(N),L)\ge 2$.
\end{enumerate}
\end{theorem}

As a corollary of Theorem \ref{alternating link}, we have the following.

\begin{corollary}\label{non-split}
Let $F$ be a closed surface embedded in $S^3$, $L$ a link contained in $F\times [-1,1]$ which has a reduced, prime, alternating diagram on $F$.
Then, $L$ is not splittable.
\end{corollary}

The proof of Theorem \ref{alternating link} is very similar to one of Theorem \ref{alternating} and the proof of Corollary \ref{non-split} is elementary, so we leave them to the readers.

\section{Preliminary}


\begin{lemma}\label{trivial}
A knot $K$ in $S^3$ is trivial if and only if any incomperssible and $\partial$-incompressible surface in the exterior $E(K)$ is an only disk.
\end{lemma}

\begin{proof} (of Lemma \ref{trivial})
Suppose that $K$ is trivial.
Then, $E(K)$ is a solid torus.
It is well-known that any incompressible and $\partial$-incompressible surface in a solid torus is an only disk.

Conversely, suppose that any incomperssible and $\partial$-incompressible surface in $E(K)$ is an only disk.
Since any knot exterior contains an incompressible non-separating orientable surface coming from a Seifert surface of minimal genus, the supposition implies that $K$ bounds a Seifert surface of genus 0.
Thus, $K$ is trivial.
\end{proof}

\begin{lemma}\label{boundary parallel}
Let $K$ be a knot in $S^3$ and $F$ an incompressible orientable surface properly embedded in $E(K)$.
If $F$ is $\partial$-compressible in $E(K)$, then $F$ is a $\partial$-parallel annulus.
\end{lemma}

\begin{proof} (of Lemma \ref{boundary parallel})
Suppose that $F$ is $\partial$-compressible in $E(K)$, and let $D$ be a $\partial$-compressing disk for $F$ in $E(K)$.
Put $\partial D\cap F=\alpha$ and $\partial D\cap \partial E(K)=\beta$.
Since $F$ is incompressible in $E(K)$, $\partial F$ consists of mutually parallel loops that are essential in the torus $\partial E(K)$.
Let $A$ be an annulus as the closure of a component of $\partial E(K)-\partial F$ containing $\beta$.

If $\beta$ is inessential in $A$, then we isotope $D$ so that $\partial D$ is entirely contained in $F$.
Then, $\partial D$ bounds a disk in $F$ since $F$ is incompressible in $E(K)$.
This shows that $\alpha$ is not essential in $F$, and contradicts that $D$ is a $\partial$-compressing disk for $F$ in $E(K)$.

Otherwise, by cutting $A$ along $\beta$ and pasting two parallel copies of $D$, we obtain a new disk $E$ with $\partial E\subset F$.
Since $F$ is incompressible in $E(K)$, $\partial E$ bounds a disk $E'$ in $F$.
Moreover, a 2-sphere $E\cup E'$ bounds a 3-ball in $E(K)$ since $E(K)$ is irreducible.
This shows that $F$ is an annulus which is $\partial$-parallel to $A$.
\end{proof}

\begin{lemma}\label{representativity}
Let $K$ be a knot contained in a closed surface $F$ in $S^3$.
Then,
\begin{enumerate}
	\item[(1)] $r(F,K)\ge1$ if and only if $F\cap E(K)$ is incompressible in $E(K)$.
	\item[(2)] $r(F,K)\ge2$ if and only if $F\cap E(K)$ is incompressible and $\partial$-incompressible in $E(K)$.
\end{enumerate}
\end{lemma}

\begin{proof} (of Lemma \ref{representativity})

(1) By the definition, $r(F,K)=0$ if and only if there exists a compressing disk $D$ for $F$ such that $|\partial D\cap K|=0$.
Equivalently, $F-K$ is compressible in $S^3-K$.

(2) Suppose that $r(F,K)=1$ and let $D$ be a compressing disk for $F$ such that $|\partial D\cap K|=1$.
Then, we have a $\partial$-compressing disk $D\cap E(K)$ for $F\cap E(K)$.
Together with (1), we have that if $F\cap E(K)$ is incompressible and $\partial$-incompressible in $E(K)$, then $r(F,K)\ge2$.

Conversely, suppose that $F\cap E(K)$ is incompressible and $\partial$-compressible in $E(K)$.
Then, by Lemma \ref{boundary parallel}, $F\cap E(K)$ is a $\partial$-parallel annulus.
This shows that $F$ is a torus and $K$ goes around $F$ once.
Hence, there exists a compressing disk $D$ for $F$ such that $|\partial D\cap K|=1$, and $r(F,K)=1$.
Thus, we have that if $r(F,K)\ge2$, then $F\cap E(K)$ is incompressible and $\partial $-incompressible in $E(K)$.
\end{proof}

\section{Proof}

\begin{proof} (of Theorem \ref{non-parallel})
Suppose that $K$ is parallel to $S$ and let $A$ be an annulus connecting $K$ and a loop in $S$.
We will show that $A$ can be isotoped so that $A\cap F=K$ and $A\cap D=\emptyset$ for a compressing disk $D$ for $F$ in $S^3-S$ on the side of containing $A$.
This implies that $F-K$ is compressible in $S^3-K$, a contradiction.

%

Let $M$ be the closure of a component of $S^3-S$ containing $K$.
Hereafter, we denote $X\cap E(K)$ by $X(K)$ for a submanifold $X$ in $S^3$.
By Lemma \ref{representativity}, we have the following claim.

\begin{claim}\label{boundary-incompressible}
$F(K)$ is incompressible and $\partial$-incompressible in $E(K)$.
\end{claim}

\begin{claim}\label{non-trivial}
$K$ is not a trivial knot in $S^3$.
\end{claim}

\begin{proof} (of Claim \ref{non-trivial})
By Claim \ref{boundary-incompressible}, $F(K)$ is an incompressible and $\partial$-incompressible surface in $E(K)$, and it is not a disk since $F$ has a positive genus.
Hence, by Lemma \ref{trivial}, $K$ is not trivial.
\end{proof}

\begin{claim}\label{annulus}
$A(K)$ is incompressible and $\partial$-incompressible in $M(K)$.
\end{claim}

\begin{proof} (of Claim \ref{annulus})
If $A(K)$ is compressible in $M(K)$, then $K$ bounds a disk in $M$.
Thus, $K$ is trivial in $S^3$ and this contradicts Claim \ref{non-trivial}.
In addition, since two boundary components of $A(K)$ are contained in the distinct surfaces $S$ and $\partial E(K)$ of $\partial M(K)$, it is $\partial$-incompressible in $M(K)$.
\end{proof}

\begin{claim}\label{A and F}
We may assume that $A\cap F$ consists of essential loops in both of $A$ and $F$.
\end{claim}

\begin{proof} (of Claim \ref{A and F})
Since both of $A(K)$ and $F(K)$ are incompressible and $\partial$-incompressible in $M(K)$, we may assume that $A(K)\cap F(K)$ consists of essential loops and arcs in both of $A(K)$ and $F(K)$.
Moreover, there exists no arc of $A(K)\cap F(K)$ which is essential in $A(K)$ since $\partial F(K)\subset \partial E(K)$ and $\partial A(K)\not\subset \partial E(K)$.
Hence, $A(K)\cap F(K)$ consists of essential loops in both of $A(K)$ and $F(K)$.
\end{proof}

\begin{claim}\label{A and F=K}
We may assume that $A\cap F=K$.
\end{claim}

\begin{proof} (of Claim \ref{A and F=K})
Suppose that $A\cap F\ne K$.
Then, by Claim \ref{A and F}, $(A\cap F)-K$ consists of essential loops in both of $A$ and $F$, and $F$ cuts $A$ into some annuli.
Let $A_1$ be the subannulus of $A$ nearest to $K$.
Since $F$ is compressible in $S^3-S$ on both sides, there exists a compressing disk $D$ for $F$ in the closure $M'$ of a component of $M-F$ containing $A_1$.
We note that by Claim \ref{annulus} and \ref{A and F}, $A_1$ is also incompressible in $M'$.
Hence, by an isotopy of $D$, we may assume that $D\cap A_1$ consists of arcs.
Moreover, by exchanging $D$ for another compressing disk if necessary, we may assume that any arc of $D\cap A_1$ is essential in $A_1$.
Then, an outermost disk on $D$ gives a $\partial$-compressing disk for $A_1$ in $M'$.
It follows that $A_1$ is parallel to a subannulus in $F$ since $F-K$ is incompressible in $M'-K$.
Thus, we can reduce $|A\cap F|$ and eventually have $A\cap F=K$.
\end{proof}

Let $D$ be a compressing disk for $F$ in $S^3-S$ on the side of containing $A$.
By Claim \ref{annulus}, we may assume that $D\cap A$ consists of arcs.
Moreover, by exchanging $D$ if necessary, we have $D\cap A=\emptyset$.
This shows that $D$ is a compressing disk for $F-K$ in $S^3-K$, and contradicts that $r(F,K)\ge2$.
\end{proof}

\begin{proof} (of Theorem \ref{alternating})

\begin{claim}\label{open disk}
$F-\pi(K)$ consists of open disks.
\end{claim}

\begin{proof} (of Claim \ref{open disk})
Suppose that there exists a component $R$ of $F-\pi(K)$ which is not an open disk.
Then, there exists an essential loop $l$ in $R$.
Let $\alpha$ be an arc in the closure of $R$ connecting $l$ and a point of $\pi(K)$ except for crossing points, and let $l'$ be a loop of $\partial N(l\cup \alpha;F)$ which intersects $\pi(K)$ in two points.
Then, $l'$ bounds a disk $D$ in $F$ such that $D\supset l\cup \alpha$ and $D\cap \pi(K)$ consists of an embedded arc since $\pi(K)$ is prime.
It follows that $l$ bounds a disk in $R$, a contradiction.
\end{proof}

We construct a graph $G$ embedded in $F$ from $\pi(K)$ in the following way.
Put a vertex of $G$ in each region of $F-\pi(K)$, and connect two vertices by an edge when the corresponding two regions face each other.

\begin{claim}\label{bipartite}
$G$ is a bipartite graph.
\end{claim}

\begin{proof} (of Claim \ref{bipartite})
We assign vertices of type (B) to a subset $X$ of the vertex set of $G$ and of type (W) to $Y$.
See Figure \ref{type} about types for a vertex.

\begin{figure}[htbp]
	\begin{center}
		\includegraphics[trim=0mm 0mm 0mm 0mm, width=.7\linewidth]{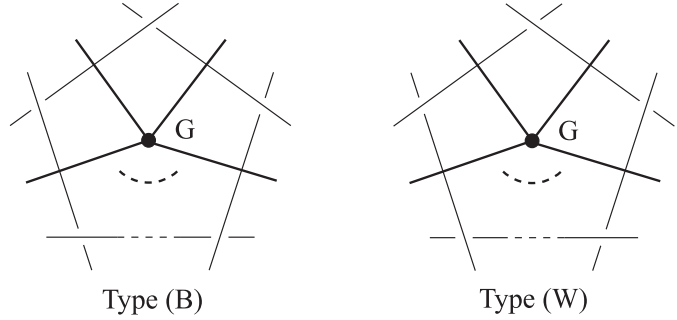}
	\end{center}
	\caption{Type for a vertex}
	\label{type}
\end{figure}

Then, we have that $X\cup Y$ is equal to the vertex set of $G$, $X\cap Y=\emptyset$ and any edge connects $X$ and $Y$ since $\pi(K)$ is alternating.
Hence, $G$ is bipartite.
\end{proof}

By Claim \ref{bipartite}, we have the conclusion 2 of Theorem \ref{alternating}.
Thus, we have two spanning surfaces $B$ and $W$ for $K$ made from the black regions and white regions.

\begin{claim}\label{non-orientable}
At least one of $B$ and $W$ is non-orientable.
\end{claim}

\begin{proof} (of Claim \ref{non-orientable})
Suppose that both of $B$ and $W$ are orientable and we assign orientations.
Then, the orientation of $B$ induces an orientation of $\pi(K)$ and each crossing has a negative sign.
Here, we note that the sign of crossings does not depend on the choice of orientations of $B$.
On the other hand, the orientation of $W$ induces an orientation of $\pi(K)$ and each crossing has a positive sign.
This is a contradiction.
Hence, at least one of $B$ and $W$ is non-orientable.
\end{proof}

By Claim \ref{non-orientable}, $K$ bounds a non-orientable surface $N$ coming from the checkerboard coloring.
Naturally, $K$ is contained in $\partial N(N)$, and $\partial N(N)-K$ is connected since $N$ is non-orientable.

\begin{claim}\label{I-bundle}
$\partial N(N)-K$ is incompressible in $N(N)$.
\end{claim}

\begin{proof} (of Claim \ref{I-bundle})
We regard the regular neighborhood $N(N)$ as a twisted $I$-bundle $N\tilde{\times}I$ over $N$.
It suffices to show that the associated $\partial I$-bundle over $N$ is incompressible in $N\tilde{\times}I$.
Consider the sequence $\pi_1(N\tilde{\times} \partial I)\to \pi_1(N\tilde{\times}I)\to \pi_1(N)$ induced by the inclusion $N\tilde{\times}\partial I\hookrightarrow N\tilde{\times}I$ and the projection $N\tilde{\times}I\to N$.
Then, the second map is an isomorphism since $N\tilde{\times}I\to N$ is a deformation retraction, and the composition is injective since $N\tilde{\times}\partial I\to N$ is a 2-fold covering.
Hence, $\pi_1(N\tilde{\times}\partial I)\to \pi_1(N\tilde{\times}I)$ is a monomorphism, and $N\tilde{\times}\partial I$ is incompressible in $N\tilde{\times}I$.
\end{proof}

In the following, we show that $\partial N(N)-K$ is incompressible in the outside of $N(N)$.
We regards $N(N)$ as the following.
For each crossing $c_i$ of $\pi(K)$, we insert a small 3-ball $B_i$ as a regular neighborhood of $c_i$.
In the rest of these 3-balls, we consider the product $R_i\times I$ for each region $R_i$ of $F-\pi(K)$.
Then, the union of $B_i$'s and $R_i\times I$'s is homeomorphic to $N(N)$.
See Figure \ref{N(N)}.

\begin{figure}[htbp]
	\begin{center}
		\includegraphics[trim=0mm 0mm 0mm 0mm, width=.6\linewidth]{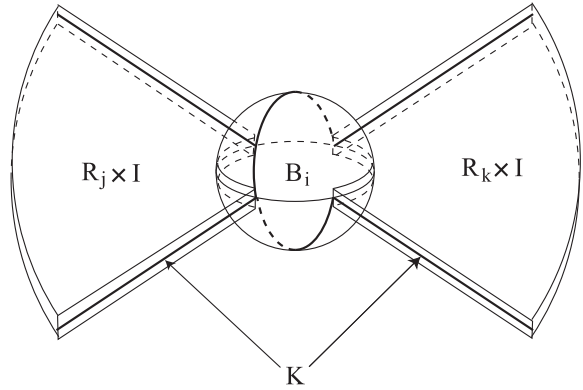}
	\end{center}
	\caption{$B_i$, $R_i\times I$ and $K$}
	\label{N(N)}
\end{figure}

Suppose that $\partial N(N)-K$ is compressible in the outside of $N(N)$ and let $D$ be a compressing disk.
Let $\Delta=\Delta_1\cup \cdots \cup \Delta_n$ be the union of components of $F-int N(N)$, where $\Delta_i$ is a disk by Claim \ref{open disk}.
Then, each component of $(\partial N(N)-K)-\partial \Delta$ is an open disk containing $R_i\times\{0\}$ or $R_i\times\{1\}$ for some $i$ and whose closure is denoted by $R_i^-$ or $R_i^+$ respectively.
Put $R=(\bigcup_i R_i^-)\cup (\bigcup_i R_i^+)$ and $B=\bigcup_i B_i$.

\begin{claim}\label{standard}
We may assume the following.
\begin{enumerate}
	\item $\partial D\cap R$ consists of arcs that connect different arc components of $\partial B\cap \partial \Delta$.
	\item $D\cap \Delta$ consists of arcs that connect different arc components of $\partial B\cap \partial \Delta$.
\end{enumerate}
\end{claim}

\begin{proof} (of Claim \ref{standard})
If $D\cap \Delta=\emptyset$, then $\partial D\subset R_i^{\pm}$ for some $i$.
This contradicts that $D$ is a compressing disk.
Therefore, $D\cap \Delta\ne\emptyset$ for any compressing disk $D$, and $\partial D\cap R$ consists of arcs.

Suppose that there exists an arc of $\partial D\cap R_i^{\pm}$ for some $i$ that connects the same arc component $\gamma$ of $\partial B\cap \partial \Delta$.
Let $\alpha$ be an outermost arc of $\partial D\cap R_i^{\pm}$ in $R_i^{\pm}$ with respect to $\gamma$, and $\delta$ the corresponding outermost disk in $R_i^{\pm}$.
Then, by an isotopy of $D$ along $\delta$, we can eliminate $\alpha$.
Hence, we may assume the condition 1 of Claim \ref{standard}.

Suppose that there exists a loop component of $D\cap \Delta$.
Let $\alpha$ be an innermost loop of $D\cap \Delta$ in $\Delta$, and $\delta$ the corresponding innermost disk in $\Delta$.
Then, by cutting $D$ along $\alpha$ and pasting $\delta$, we have a new compressing disk and a 2-sphere, and $\alpha$ is eliminated.
Hence, we may assume that $D\cap \Delta$ consists of arcs.

Next, suppose that there exists an arc of $D\cap \Delta$ that connects the same arc component $\gamma$ of $\partial B\cap \partial \Delta$.
Let $\alpha$ be an outermost arc of $D\cap \Delta$ in $\Delta$ with respect to $\gamma$, and $\delta$ the corresponding outermost disk in $\Delta$.
Then, by cutting $D$ along $\alpha$ and pasting two parallel copies of $\delta$, we have two disks whose boundaries are contained in $\partial N(N)-K$.
Since $\partial D$ is essential in $\partial N(N)-K$, at least one of boundaries of these two disks is also essential in $\partial N(N)-K$.
Thus, we have a new compressing disk and $\alpha$ is eliminated.
Hence, we may assume the condition 2 of Claim \ref{standard}.
\end{proof}

Next, we concentrate on an outermost arc $\alpha$ of $D\cap \Delta$ in $D$ and the corresponding outermost disk $\delta$ in $D$.
Put $\delta\cap \partial D=\beta$.
By Claim \ref{standard}, we have two configurations.

\begin{description}
	\item[Case 1] $\beta$ connects the same crossing ball $B_i$ (Figure \ref{alternate1}).
	\item[Case 2] $\beta$ connects different crossing balls $B_i$ and $B_j$ (Figure \ref{alternate2}).
\end{description}

\begin{figure}[htbp]
	\begin{center}
		\includegraphics[trim=0mm 0mm 0mm 0mm, width=.7\linewidth]{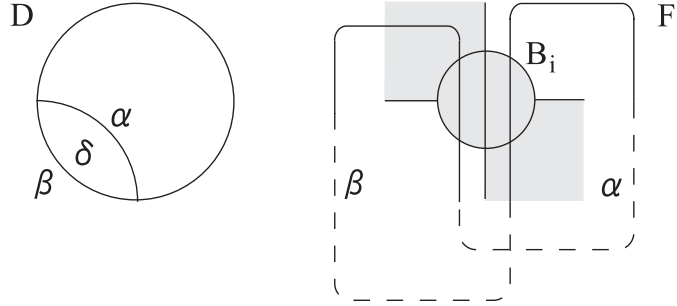}
	\end{center}
	\caption{Configuration of Case 1}
	\label{alternate1}
\end{figure}

\begin{figure}[htbp]
	\begin{center}
		\includegraphics[trim=0mm 0mm 0mm 0mm, width=.7\linewidth]{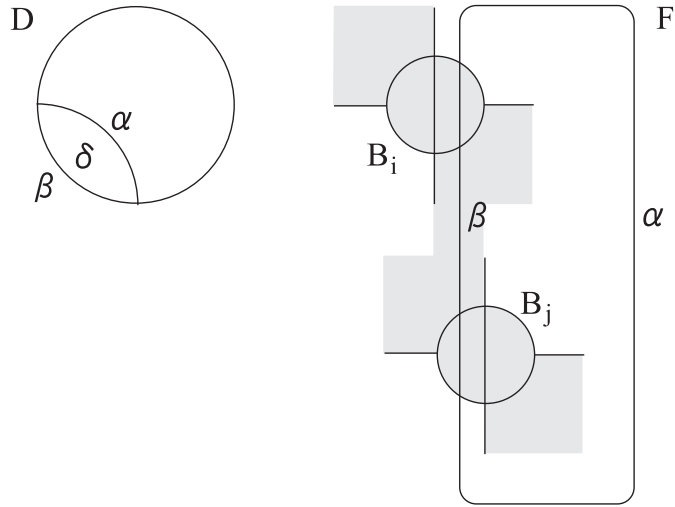}
	\end{center}
	\caption{Configuration of Case 2}
	\label{alternate2}
\end{figure}

In Case 1, by connecting $\partial \beta$ on $\partial B_i$ and projecting on $F$, we have a loop $l_\beta$ on $F$ which intersects $\pi(K)$ in one crossing point $c_i$.
Similarly, we obtain a loop $l_\alpha$ on $F$ which intersects $\pi(K)$ in one crossing point $c_i$.
Since $l_\beta$ intersects $l_\alpha$ in one point $c_i$, $l_\beta$ is essential in $F$.
Let $l_\beta$ avoid $c_i$.
Then we have an essential loop in $F$ which intersects $\pi(K)$ in two points except for crossings.
This contradicts that $\pi(K)$ is prime.

In Case 2, we have a loop $\pi(\alpha\cup\beta)$ in $F$ which intersects $\pi(K)$ in two points except for crossings.
It does not bound a disk $D$ in $F$ such that $D\cap \pi(K)$ is an arc since there are crossings $c_i$ and $c_j$ on both sides of the loop.
This contradicts that $\pi(K)$ is prime.

Hence, $\partial N(N)-K$ is incompressible in the outside of $N(N)$, and together with Claim \ref{I-bundle}, $\partial N(N)-K$ is incompressible in $S^3-K$.
By Lemma 2.3 (1), we have $r(\partial N(N),K)\ge1$.

Next, we show that $r(\partial N(N),K)\ge 2$.
Suppose that $r(\partial N(N),K)=1$, equivalently, that $\partial N(N)\cap E(K)$ is $\partial$-compressible in $E(K)$.
By Lemma \ref{boundary parallel}, $\partial N(N)\cap E(K)$ is a $\partial$-parallel annulus.
Therefore, $\partial N(N)$ is a torus, $N(N)$ is a solid torus, and $N$ is a M\"{o}bius band.
Moreover, since $\pi(K)$ is reducible, $\pi(K)$ is a standard $(2,n)$-torus knot diagram, where $n$ is an odd integer.
If $|n|\ge3$, then $r(\partial N(N),K)=2$.
This contradicts the supposition $r(\partial N(N),K)=1$.
Otherwise, $n=\pm1$.
This shows that $\pi(K)$ is reducible, a contradiction.
\end{proof}

\section{Conclusion}

After Corollary \ref{alternating2}, it remains to consider whether a generalized alternating knot on a closed surface is trivial for itself.

\begin{conjecture}\label{conjecture}
Let $F$ be a closed surface embedded in $S^3$, $K$ a knot contained in $F\times [-1,1]$ which has a reduced, prime, alternating diagram on $F$.
Then, $K$ is not trivial for $F$.
\end{conjecture}

Conjecture \ref{conjecture} is true for a 2-sphere since an ordinary alternating knot is non-trivial, and for an unknotted torus since a torus knot exterior contains only two essential surfaces, the cabling annulus and the minial genus Seifert surface.

\bibliographystyle{amsplain}

\end{document}